\documentclass{amsart}

\usepackage{amsmath,amsthm,amssymb}
\usepackage{amsfonts}
\usepackage{rotating}
\usepackage{euscript}
\usepackage{pst-node}
\usepackage{epsfig}

\def\be{\begin{equation}}
\def\ee{\end{equation}}

\def\C{{\mathbb C}} 
 
\def\P{{\mathbb P}}
\def\Z{{\mathbb Z}}

\def\f{\EuScript}

\def\phi{{\varphi}}
\def\v{{\varepsilon}} 
\def\deg{{\rm deg\,}}

\def\cos{{\rm cos\,}} 
 
\def\GCD{{\rm GCD }}
\def\LCM{{\rm LCM }}

\def\qed{$\ \ \Box$ \vskip 0.2cm}

\def\bp{\begin{proposition}}
\def\ep{\end{proposition}}

\def\bt{\begin{theorem}}
\def\et{\end{theorem}}
\def\br{\begin{remark}}
\def\er{\end{remark}}
\def\be{\begin{equation}}
\def\bee{\begin{equation*}}
\def\la{\label}
\def\l{\label}

\def\ee{\end{equation}}
\def\eee{\end{equation*}}
\def\bl{\begin{lemma}}
\def\el{\end{lemma}}
\def\br{\begin{remark}}
\def\er{\end{remark}}
\def\bc{\begin{corollary}}
\def\ec{\end{corollary}}
\def\pr{\noindent{\it Proof. }}

\def\bd{\begin{definition}}
\def\ed{\end{definition}}
\input epsf.sty

\newtheorem{theorem}{Theorem}[section]
\newtheorem{lemma}[theorem]{Lemma}
\newtheorem{corollary}[theorem]{Corollary}
\newtheorem{proposition}[theorem]{Proposition}
\newtheorem{remark}{Remark}[section]

\begin{document}
\title[Polynomial moment problem]{Generalized ``second Ritt theorem''
and explicit solution of the polynomial moment problem.
}
\author{F. Pakovich}
\address{Department of Mathematics, Ben-Gurion University
of the Negev, P.O.B. 653, Beer-Sheva, Israel} 
\email{pakovich@cs.bgu.ac.il}
\date{}

\begin{abstract} In the recent paper \cite{pm} was shown that any 
solution of ``the 
polynomial moment problem'', 
which asks to describe polynomials $P,Q$ 
satisfying $\int_{a}^bP^kd Q=0$ for all $k\geq 0$, 
may be obtained as a sum of some ``reducible" solutions 
related to different decompositions 
of $P$ into a composition of two 
polynomials of lesser degrees.
However, the methods of \cite{pm} do not permit to estimate the number of necessary reducible solutions or to describe them explicitly.
In this 
paper we provide a description of the 
polynomial solutions of the functional equation 
$P=P_1\circ W_1=P_2\circ W_2=\dots =P_r\circ W_r,$
and on this base describe solutions of the polynomial moment problem in an explicit form suitable 
for applications.  

\end{abstract}

\maketitle

\section{Introduction}

About a decade ago,  
in the series of papers \cite{bfy1}--\cite{bfy4} the following
``polynomial moment problem'' was posed: for a given complex polynomial $P$ 
and complex numbers 
$a,b$
describe polynomials $Q$ satisfying the system of equations
\be\la{111}
\int^b_a P^k\,dQ \,=\, 0, \ \  \ k\geq 0.
\ee 
Despite its rather classical and simple setting this problem turned out to 
be quite difficult and was intensively studied in many recent papers 
(see, e.\,g., \cite{bfy2}--\cite{bfy5},
\cite{c}, \cite{pa1}--\cite{pp}, \cite{pm}, \cite{ro}).

The main motivation for the study of the polynomial moment problem is its relation with the
center problem for the Abel differential equation 
\be \l{ab}
\frac{d y}{d z}=p(z)y^2+q(z)y^3
\ee
with polynomial
coefficients $p, q$ in the complex domain.
For given $a,b\in \C$ the center problem for the Abel equation
is to find necessary and sufficient conditions on $p,q$ which imply
the equality $y(b)=y(a)$ for any solution $y(z)$ of \eqref{ab} with $y(a)$ small enough.
This problem is closely related to the classical 
Center-Focus problem of Poincar\'e and
has been studied in many recent papers 
(see e.g. \cite{bby}-\cite{c}, \cite{y1}).

The center problem for the Abel equation is connected with the 
polynomial moment problem in several ways.
For example, it was shown in \cite{bfy3}
that for the parametric version  
$$
\frac{d y}{d z}=p(z)y^2+\varepsilon q(z)y^3
$$
of \eqref{ab} the ``infinitesimal'' center conditions with respect to $\varepsilon$ 
reduce to equations \eqref{111}, where 
\be \la{intt} P(z)=\int p(z) d z,\ \ \ Q(z)=\int q(z) d z.\ee
On the other hand, it was shown in \cite{bry} that ``at infinity'' (under an appropriate projectivization of the 
parameter space) the system of equations on 
coefficients of $p$ and $q$ describing the center set of \eqref{ab} also reduces to equations \eqref{111}, where $P$ and $Q$ are defined as above.
Many other results relating
the center problem and the polynomial moment problem may be found in \cite{bry}.

There exists a natural condition on $P$ and $Q$
which reduces equations \eqref{111}, \eqref{ab} to similar equations
with respect to polynomials of smaller degrees. Namely, suppose that
there exist polynomials $\widetilde P,$ $\widetilde Q,$ $W,$ $\deg W>1$, such that
\be \l{2}
P=\widetilde P\circ W, \ \ \ \ \ 
\ Q=\widetilde Q\circ W,
\ee
where the symbol $\circ$ denotes a superposition of functions: 
$f_1\circ f_2=f_1(f_2)$. 
Then after the change of variable $z\rightarrow W(z)$
equations \eqref{111} transform to the equations  
\be \l{zam}
\int^{W(b)}_{W(a)} \widetilde P^k d \widetilde Q=0, \ \ \ k\geq 0, 
\ee
while equation \eqref{ab} transforms to the equation 
\be \l{ab1}
\frac{d\widetilde y}{d w}=\widetilde P^{\prime} (w)\widetilde y^2+\widetilde Q^{\prime}(w)\widetilde y^3.
\ee

Furthermore, if the polynomial $W$ in \eqref{2} 
satisfies the equality 
\be \l{w} W(a)=W(b), 
\ee 
then it follows from the Cauchy theorem that all integrals in \eqref{zam} vanish implying that all integrals in 
\eqref{111} also vanish.
Similarly, since any solution $y(z)$ of equation \eqref{ab} is the pull-back 
$$y(z)=\widetilde y(W(z))$$
of a solution $\widetilde y(w)$ 
of equation \eqref{ab1}, if $W$ satisfies \eqref{w}, then equation
\eqref{ab} has a center. 
We will call a center for equation
\eqref{ab} or a solution of system \eqref{111} 
{\it reducible} if there exist polynomials 
$\widetilde
P,$ $\widetilde Q,$ $W$ such that
conditions 
\eqref{2}, \eqref{w} hold.
The main conjecture concerning the center problem for the Abel equation
(``the composition conjecture for the Abel equation")
states 
that any center for the Abel equation is reducible 
(see \cite{bry}
and the bibliography therein).

By analogy with the composition conjecture for the Abel equation it was suggested 
(``the composition conjecture for the polynomial moment problem") 
that
any solution of (1) is reducible.
This conjecture was shown to be 
true in many cases. For instance, 
if $a,b$ are not critical points of $P$   
(\cite{c}), if $P$ is indecomposable that is can not be represented as a composition of two polynomials of lesser degrees (\cite{pa2}), and in some other 
special cases (see e. g. \cite{bfy3}, \cite{pp}, \cite{pry}, \cite{ro}).
Nevertheless, 
in general the composition conjecture for the polynomial moment problem 
fails to be true. 

A class of counterexamples to the composition conjecture for the polynomial moment problem was constructed in \cite{pa1}. These counterexamples use
polynomials $P$ which admit 
``double decompositions'' of the form
\be \la{rii} P=P_1\circ W_1=P_2\circ W_2,\ee 
where $P_1,$ $P_2,$ $W_1,$ $W_2$
are non-linear polynomials. If $P$ is such a polynomial and,
in addition, the equalities $$W_1(a)=W_1(b),\ \ \ \ W_2(a)=W_2(b)$$ hold, then
for any polynomials $V_1,V_2$ the polynomial
$$Q=V_1\circ W_1+V_2\circ W_2$$ satisfies \eqref{111} by linearity.
On the other hand, it can be shown (see \cite{pa1}) that if $\deg W_1$ and 
$\deg W_2$ are coprime, then condition \eqref{2} is not satisfied already for $Q=W_1+W_2.$
 
Notice that the description of polynomial solutions of \eqref{rii}
may be reduced to the case where
\be \la{gcd} \GCD(\deg P_1,\deg P_2)=1, \ \  \ \GCD(\deg W_1,\deg W_2)=1\ee 
(see Section 2 below). On the other hand, in the last case solutions are  described explicitly by the so called ``second Ritt theorem'' which states that if 
if $P_1,P_2,W_1,W_2$ are polynomials satisfying \eqref{rii} and \eqref{gcd}, then there exist polynomials 
$\nu,$ $\mu,$ $\sigma_1,$ $\sigma_2$ of degree one such that 
up to a possible replacement of $P_1$ to $P_2$ and $W_1$ to $W_2$ either
\begin{align*} &P_1=\nu \circ z^n \circ \sigma_1^{-1}, & &W_1=\sigma_1 \circ
z^sR(z^n) \circ \mu \\ 
&P_2=\nu \circ z^sR^n(z) \circ \sigma_2^{-1},& 
&W_2=\sigma_2 \circ z^n \circ \mu, \end{align*} 
where $R$ is a polynomial, or
\begin{align*} &P_1=\nu \circ T_n \circ \sigma_1^{-1}, &
&W_1=\sigma_1 \circ T_m \circ \mu, \\
&P_2=\nu \circ T_m \circ \sigma_2^{-1}, & &W_2=\sigma_2 \circ
T_n\circ \mu, \end{align*} 
where $T_n, T_m$ are the Chebyshev polynomials.

It was conjectured in \cite{pa4} that actually any solution 
of~\eqref{111} can be represented as a {\it sum} of reducible ones and recently
this conjecture was proved in~\cite{pm}. In more details, it was proved in~\cite{pm}
that  
non-zero polynomials $P,$ $Q$ satisfy system (1) if and only if 
\be \la{aza} Q=\sum_{i=1}^rQ_i\ee where $Q_i,$ $1\leq i \leq r,$ are polynomials 
such that
\be \l{cc}
P=P_i\circ W_i, \ \ \
Q_i=V_i\circ W_i, \ \ \  \ \ \ W_i(a)=W_i(b)
\ee
for some polynomials $P_i, V_i, W_i,$ $1\leq i \leq r.$
Although this result in a sense solves the problem it does not provide any 
explicit description of polynomials $P$ and $Q$ satisfying \eqref{aza}, \eqref{cc}.
On the other hand, for applications to differential equations (for example, for the investigation of the center set for 
the Abel equation) such a description would be highly desirable. 

The problem of explicit description of solutions of the polynomial moment problem
naturally leads to the following two problems.

First, since the number $r$ in \eqref{aza} is a priory unbounded, it is necessary to describe somehow polynomial solutions of the equation \be \la{serr} P_1\circ W_1=P_2\circ W_2=\dots =P_r\circ W_r\ee for $r>2$. Note that as in the case $r=2$ 
such a description 
reduces to the case where
\be\la{usl2} \hskip -0.8cm \GCD(\deg P_1,\deg P_2,\dots ,\deg P_r)=1, \  \GCD(\deg W_1,\deg W_2,\dots ,\deg W_r)=1.\ee  However, since 
\eqref{usl2} does not imply that 
the degrees of polynomials $P_i,$ $1\leq i \leq r,$ as well as of  
$Q_i,$ $1\leq i \leq r,$ are necessary {\it pairwise} coprime, the Ritt theorem does not provide any immediate information about solutions of \eqref{serr}. 

Second, since the representation of a solution $Q$ in the form of a sum of reducible solutions is not unique, it is desirable to find some canonical form of such 
a representation. We may illustrate it by the following example.
Let $P=T_n$. Then for any divisor $d$ of $n$ the equality $$T_n=T_{n/d}\circ T_d$$ holds and therefore any $W_i=T_{d_i}$,
where $d_i\mid n,$ is  
a compositional right 
factor of $P.$ However, one can show (see Section \ref{pp2a}) that 
if the equalities 
$$T_{d_1}(a)=T_{d_1}(b), \ \ \ \ \ \ T_{d_2}(a)=T_{d_2}(b), \ \ \ \ \ \ 
T_{d_3}(a)=T_{d_3}(b)$$ hold, then there exists a pair of indices 
$d_{i_1},d_{i_2},$ $i_1\neq i_2,$ such that the polynomials $T_{d_{i_1}}$ and $T_{d_{i_2}}$ have a common 
compositional right factor $T_l$ such that $T_l(a)=T_l(b).$ Therefore, in 
any solution $$Q=\sum_{i=1}^rQ_i=\sum_{i=1}^rV_{i}\circ T_{d_{i}}$$ such that $r>2$ 
we may replace the sum of two reducible solutions  
$$Q_{i_1}+Q_{i_2}=V_{i_1}\circ T_{d_{i_1}}+V_{i_2}\circ T_{d_{i_2}}$$ by the unique reducible solution 
$$(V_{i_1}\circ T_{d_{i_1}/l}+V_{i_2}\circ T_{d_{i_1}/l})\circ  T_{l},$$
and continuing this process eventually represent $Q$ as a sum of at most {\it two} reducible solutions.

In this paper we solve both problems above. Our first result is an analogue of 
the second Ritt theorem for functional equation \eqref{serr}. 
Recall that two polynomials $U$, $V$ are called linearly equivalent if $U=\mu\circ V\circ\nu$ for some polynomials $\mu,\nu$ of degree one.

\bt \la{os1} Let $P_i,W_i,$ $1\leq i \leq r,$ be polynomials satisfying \eqref{serr} and
\eqref{usl2}. Then at least one $P_i,$ $1\leq i \leq r,$ 
is linearly equivalent either to a Chebyshev polynomial
or to a power, and at least one $W_i,$ $1\leq i \leq r,$ is linearly equivalent either to a Chebyshev polynomial
or to a power.
\et

Note that although in distinction with the second Ritt 
theorem this result does not provide a full description of {\it all} polynomials involved in \eqref{serr}, it still implies their ``partial'' description sufficient 
for applications (see Subsection \ref{appp}).

Theorem \ref{os1} permits to reduce the number of reducible solutions in the 
representation $Q=\sum_{i=1}^rQ_i$ in the way similar to the one described above and eventually to show that any non-reducible solution of the polynomial moment problem
may be represented either as a sum of two reducible solutions 
related to double decompositions
appearing in the second Ritt theorem
or as a sum of three reducible solutions
related to a special ``triple'' decomposition which may be described as
follows. 

Let
$$P=z^{2}R^2(z^{2})\circ T_{mn},$$ where $R$ is a polynomial and $m,n$ are odd numbers such that $\GCD(m,n)=1.$  
Then $W_1=T_{2n}$ and $W_2=T_{2m}$ are compositional right factors of
$P$ since $$P=zR^{2}(z)\circ z^2\circ T_{mn}=zR^{2}(z)\circ \frac{z+1}{2}\circ T_2
\circ T_{mn}=\frac{z+1}{2}R^{2}\left(\frac{z+1}{2}\right)\circ T_{2mn}=$$ 
$$=\frac{z+1}{2}R^{2}\left(\frac{z+1}{2}\right)\circ  T_{m}\circ T_{2n}=
\frac{z+1}{2}R^{2}\left(\frac{z+1}{2}\right)\circ  T_{n}\circ T_{2m}.$$
Furthermore, since 
$$P=z^{2}\circ zR(z^{2})\circ T_{mn},$$ the polynomial
$W_3=zR(z^2) \circ T_{mn}$ also is a compositional right factor
of $P$, and one can show that induced solutions of the polynomial moment problem in general can not be reduced to solutions related to the Ritt theorem.

More precisely, our principal result concerning the polynomial moment problem is the following theorem.

\bt \l{os2} Let $P,$ $Q$ be non-constant polynomials and $a,b$ be distinct complex numbers satisfying system (1).
Then one of the following conditions holds:
\vskip 0.25cm 
\noindent 1) There exist polynomials $\widetilde P,$ $\widetilde Q,$ $W$
such that 
$$ P=\widetilde P\circ W, \ \ \ \ Q=\widetilde Q\circ W, \ \ \ {\it and} \ \ \  
W(a)=W(b);$$ 
\vskip 0.25cm 
\noindent 2) There exist polynomials $V_1,$ $V_2,$ $R,$ $W,$ $U$ such that 
$$P=U\circ z^{sn}R^n(z^n) \circ W,\ \ \ Q=V_1\circ z^n \circ W+ 
V_2\circ  z^sR(z^n)\circ W,$$ and 
$$W^n(a)=W^n(b), \ \ \ R(W^n(a))=0,$$ 
where $n>1$, $s\geq 1,$ $\GCD(s,n)=1;$
\vskip 0.25cm 
\noindent 3) There exist polynomials $V_1,$ $V_2,$ $U,$ $W$ and the Chebyshev polynomials $T_n,$ $T_m$, $T_{nm}$ such that 
$$P=U\circ T_{nm} \circ W,\ \ \ Q=V_1\circ T_n \circ W+ 
V_2\circ  T_m\circ W,$$ and 
$$T_n(W(a))=T_n(W(b)), \ \ \ T_m(W(a))=T_m(W(b)),$$
where $n>1,$ $m>1$, $\GCD(m,n)=1;$
\vskip 0.25cm 
\noindent 4) There exist polynomials $V_1,$ $V_2,$ $V_3,$ $U,$ $W,$ $R$ and the Chebyshev polynomials $T_n,$ $T_m$, $T_{nm}$ such that 
$$P=U\circ z^{2}R^2(z^{2})\circ  T_{mn}
\circ W,$$ $$Q=V_1\circ T_{2n} \circ W+ 
V_2\circ  T_{2m}\circ W+V_3\circ (zR(z^2) \circ T_{mn}) \circ W,$$ and 
$$T_{n}(W(a))=-T_{n}(W(b)), \ \ \ T_{m}(W(a))=-T_{m}(W(b)),$$ $$W(a)\neq - W(b),\ \ \  R(1)=0,$$
where $n>1,$ $m>1$ are odd and $\GCD(m,n)=1.$

\et

Notice that the requirements imposed on $a,$ $b$ in case 4) obviously imply the equalities $$T_{2n}(W(a))=T_{2n}(W(b)),\  \ \  T_{2m}(W(a))=T_{2m}(W(b)).$$
Furthermore, one can show (see Section 4) that these requirements imply the equality 
$T_{mn}^2(a)=T_{mn}^2(b)=1$ and hence the equality  
$$(zR(z^2) \circ T_{mn})(a)=(zR(z^2) \circ T_{mn})(b).$$
Therefore, as in other cases, $Q$ is 
exactly a sum of reducible solutions. The additional restrictions are imposed since otherwise 4) reduces to 2) or 3).

The paper is organized as follows. In the second section we recall 
the description of polynomial solutions of equation \eqref{rii}.
In the third section we prove Theorem \ref{os1}.
Finally, in the fourth section we 
prove Theorem \ref{os2}.

\section{Polynomial solutions of $P_1\circ W_1=P_2\circ W_2$}

\subsection{Imprimitivity systems and decompositions of rational functions}
In this subsection we recall the correspondence between equivalence classes of  decompositions of a rational function $F$ into compositions of rational functions of lesser degrees and imprimitivity systems of the monodromy group of $F$. For a more detailed account of algebraic structures related to decompositions of rational functions see e.g. \cite{mp}, Section 2.1. 

Let $G\subseteq S_n$ be a transitive permutation group acting on the set 
$X=\{1,2,\dots ,n\}$.
A proper subset $B$ of $X$ 
is called a {\it block} of $G$ if for each $g\in G$ the set
$B^g$ is either disjoint or equal to $B$ (see e.g. \cite{wi}). For a block $B$ the sets
$B^g,$ $g\in G,$ form a partition
of $X$ into a disjoint union
of blocks of equal cardinality which is called an {\it imprimitivity system} of $G$.

If $F$ is a rational function with complex coefficients, then the structure of decompositions of $F$ into compositions of rational functions of lesser degrees
is defined by the structure of imprimitivity systems of its monodromy group $G.$  
Namely, suppose that 
$G$ is realized as a permutation group acting
on the set $F^{-1}\{z_0\}$, consisting of preimages of a non critical value $z_0$ of $F$
under the map $F\,:\C\P^1\rightarrow \C\P^1$. Further, let $F=A\circ B$ be a decomposition of $F$ and $x_1,x_2,\dots,x_r$ be preimages of $z_0$ under the map $A\,:\C\P^1\rightarrow \C\P^1$. Then the sets
$X_i=B^{-1}\{x_i\},$ $1\leq i\leq r,$  
form an imprimitivity system $\f E$ of $G$. 
Furthermore, if $\f E$ and $\widetilde{\f E}$ are imprimitivity systems corresponding to 
decompositions $A\circ B$ and $\widetilde A\circ \widetilde B$ of $F$, then  
$\f E$ is a refinement of $\widetilde{\f E}$ if and only if there exists a rational function $U$
such that $$A=\widetilde A\circ U, \ \ \ U\circ B=\widetilde{ B}.$$ In particular, $\f E=\widetilde{\f E}$ if and only if there exists a rational function of degree one $\mu$ such that 
\be \la{verb} \widetilde A=A\circ \mu, \ \ \ \ B=\mu^{-1}\circ \widetilde B.\ee
In the last case we will call decompositions $F=A\circ B$ and $F=\widetilde A\circ \widetilde B$ {\it equivalent} and will use the notation $A\circ B\sim \widetilde A\circ \widetilde B.$

It is easy to see that if $F=A\circ B$ is a decomposition of a {\it polynomial} into a composition of rational functions, then $A\circ B\sim \widetilde A\circ \widetilde B,$ where $\widetilde A,$ $\widetilde B$ are polynomials.
Taking into account this fact, below we always will assume that all the functions
considered are polynomials and will use the following modification of the general definition of equivalence: two decompositions of a {\it polynomial} $F$ into compositions of {\it polynomials} $F=A\circ B$ and $F=\widetilde A\circ \widetilde B$ are called equivalent if there exists a {\it polynomial} of degree one $\mu$ such that \eqref{verb} holds.

\subsection{\la{iiii} Chebyshev polynomials and their properties}
Let $U$, $V$ be polynomials. We will say that $U$ is linearly equivalent to $V$ and
will use the notation $U\sim V$ if $U=\mu\circ V\circ\nu$ for some polynomials $\mu,\nu$ of degree one.
In this subsection we recall the definition of Chebyshev polynomials and their 
characterization up to the linear equivalence.

Let $P$ be a polynomial of degree $n$ and 
$S(P)=\{z_1,z_2, \dots ,z_s\}$ be the ordered set 
of all {\it finite} 
critical values of $P$. Clearly, for each $j,$ $1\leq j \leq s,$ 
the set $$\Pi_j(P)=(a_{j,\,1},a_{j,\,2}, ... , a_{j,\,p_j}),$$ consisting of local multiplicities of $P$ at the points of $P^{-1}\{z_j\},$
is a partition of the number $n$. The collection of partitions 
$$\Pi(P)=\{\Pi_1(P), \Pi_2(P), \dots \Pi_s(P)\}$$ 
is called the {\it passport} of $P$.
It follows from the Riemann existence theorem that a polynomial $P$ is defined by the sets $S(P)$ and $\Pi(P)$ up to 
the change $P\rightarrow P\circ \mu$, where $\mu$ is a polynomial of degree one.
Note that this implies in particular that if $S(P)$ contains only two points, then $P$ is defined by its   
passport up to the linear equivalence. Note also that if $S(P)$ contains only one point, then $P\sim z^n$ for some $n \geq 1.$

The Chebyshev polynomials may be defined by 
the formula \be \la{cheee} T_n(\cos \phi)=\cos n\phi.\ee
It follows easily from this definition that 
if $n>2$, then $S(T_n)=\{-1,1\}$ 
and $\Pi(T_n)$ is \be \la{pas1} \{(1,1,2,2,\dots 2),\ (2,2,\dots ,2)\},\ee if $n$ is even, or 
\be \la{pas2} \{(1,2,2,\dots 2),  \ (1,2,2,\dots ,2)\},\ee if $n$ is odd. 
Furthermore, in view of the remark above, a Chebyshev polynomial is defined by 
its passport up to the linear equivalence.
In particular, $T_n$ is not linearly equivalent to $z^n$ unless $n=2.$

Finally, notice that \eqref{cheee} implies 
the equality \be \la{ioi} T_n(-z)=(-1)^nT_n(z),\ \ \ n\geq 1.\ee

\subsection{\la{rittt} Reduction to the case of coprime degrees. }
The description of polynomial solutions of the equation 
\be \la{ura} P_1\circ W_1=P_2\circ W_2\ee 
may be reduced to the case where 
\be \la{ass0} \GCD(\deg P_1,\deg P_2)=1,  \ \ \ \GCD(\deg W_1,\deg W_2)=1\ee
owing to the statement given below. Since in the following we will need a generalization 
of this statement, we provide its complete proof.

\bt {\rm (}\cite{en}, \cite{tor}{\rm )}\la{r1}
Let $P_1,P_2,W_1,W_2$ be polynomials 
such that \eqref{ura} holds. Then there exist polynomials
$U, V, \widetilde P_1, \widetilde P_2, \widetilde W_1, \widetilde W_2,$ where
$$\deg U=\GCD(\deg P_1,\deg P_2),  \ \ \ \deg V=\GCD(\deg W_1,\deg W_2),$$
such that
$$P_1=U\circ \widetilde P_1, \ \  P_2=U\circ \widetilde P_2, \ \ W_1=\widetilde W_1\circ V, \ \  W_2=\widetilde W_2\circ V,$$
and 
$$ \widetilde P_1\circ \widetilde W_1=\widetilde P_2\circ \widetilde W_2.$$\et
\vskip 0.2cm
\noindent{\it Proof.}
Let $P$ be the polynomial defined by equality \eqref{ura}.
Clearly, the monodromy group $G$ of $P$ contains a cycle $\sigma$ of length $n=\deg P$, corresponding to a loop 
around infinity, and without loss of generality we may assume that this cycle coincides with the cycle 
$\sigma=(12...n)$. 

Since any $\sigma$-invariant partition of the set
$\{1,2,\dots,n\}$ coincides with the set $I_{d}$
consisting of residue classes 
modulo $d$ for some $d\vert n$,
any imprimitivity system of $G$ also has such a form. Owing to the correspondence between decompositions of $P$ and imprimitivity systems of $G$ this implies easily that 
Theorem \ref{r1} is equivalent to the following statement: if $I_{d_1}$ and
$ I_{d_2}$ are imprimitivity systems of $G$ corresponding to divisors 
$d_1$ and $d_2$ of $n$ respectively, then $I_{\LCM(d_1,d_2)}$ and $I_{\GCD(d_1,d_2)}$  
also are imprimitivity systems of $G.$ 

In order to prove the first part of the last statement observe that for any element $x\in X$
the intersection of two blocks $B_1$, $B_2$ containing $x$ obviously is a block and,  
if $B_1\in I_{d_1},$ $B_2\in I_{d_2},$ then $B_1\cap B_2$ coincides
with a residue class modulo $\LCM(d_1,d_2).$
The easiest way to prove the second part is to observe that $I_{d}$ is 
an imprimitivity system for $G$ if and only if the subspace $V_d$ of $\C^{n}$,
consisting of $d$-periodic vectors, is invariant with respect to the 
permutation representation $\rho_G$
of $G$ on $\C^n$, where by definition for $g\in G$ and 
$\vec v \in \C^n,$ $\vec v=(a_1,a_2,\dots ,a_n),$ the vector  
$\vec v^g$ is defined by the formula $\vec v^g=(a_{1^g},a_{2^g},\dots ,a_{n^g})$ (see \cite{pm}, Section 3.1). Clearly, if $V_{d_1}$ and 
$V_{d_2}$ are $\rho_G$-invariant, then the subspace 
$V_{d_1}\cap V_{d_2}$ also is $\rho_G$-invariant. 
On the other hand, it is easy to see that $V_{d_1}\cap V_{d_2}=V_{\GCD(d_1,d_2)}.$  \qed

Let us mention the following well known corollaries of Theorem \ref{r1}.
\bc \la{c0} Let $P_1,P_2,W_1,W_2$ be polynomials 
such that \eqref{ura} holds. Assume additionally that $\deg W_1\vert\deg W_2$ or equivalently that 
$\deg P_2\vert\deg P_1$. Then 
there exists a polynomial $S$ such that 
$$P_1=P_2\circ S, \ \ \ W_2=S\circ W_1.$$ In particular, if $\deg W_1=\deg W_2,$ then 
there exists a polynomial $\mu$ of degree one such that 
$$P_1=P_2\circ \mu, \ \ \ W_1=\mu^{-1}\circ W_2.$$
\ec
\pr Indeed, if $\deg W_1\vert\deg W_2$, then the degree of the polynomial $\widetilde W_1$ from 
Theorem \ref{r1} is one and hence the equality $W_2=S\circ W_1$ holds for $S=\widetilde W_2\circ \widetilde W_1^{-1}$.
Now the equality 
$$P_1\circ W_1=P_2\circ W_2=P_2\circ S\circ W_1$$ implies that $P_1=P_2\circ S.$
\qed

\bc \la{c1}
Let $P_1,W_1$ be polynomials such that $P_1\circ W_1=z^n.$ Then there exists a polynomial $\mu$ of degree one such that 
$$P_1=z^d\circ \mu, \ \  \ W_1=\mu^{-1}\circ z^{n/d}$$ for some $d\vert n.$ Similarly, if $P_1\circ W_1=T_n,$ then there exists a polynomial $\mu$ of degree one such that 
$$P_1=T_d\circ \mu, \ \  \ W_1=\mu^{-1}\circ T_{n/d}$$ for some $d\vert n.$  
\ec
\pr Clearly, any of the equalities $P_1\circ W_1=z^n$ and $P_1\circ W_1=T_n$ implies 
that $d=\deg P_1$ is a divisor of $n.$ 
On the other hand, for any $d\vert n,$ the equalities
$$z^n=z^d\circ z^{n/d}, \ \ \ T_n=T_d\circ T_{n/d}$$ hold.
Therefore, Corollary \ref{c1} follows from Corollary \ref{c0} applied 
to the equalities $P_1\circ W_1=T_d\circ T_{n/d}$ and $P_1\circ W_1=z^d\circ z^{n/d}$. \qed

\subsection{The second Ritt theorem}
An explicit description of polynomials satisfying \eqref{ura}, \eqref{ass0} is given by the following statement known as the second Ritt theorem (see \cite{r1} 
as well as \cite{bilu}, \cite{f1}, \cite{plau}, \cite{sch}, \cite{tor}, \cite{za}, \cite{mz}).
 
\bt {\rm (}\cite{r1}{\rm )} \la{r2}
Let $P_1,P_2,W_1,W_2$ be polynomials 
such that \eqref{ura} and \eqref{ass0} hold.
Then there exist polynomials $\sigma_1,\sigma_2,\mu, \nu$ of degree one  
such that, up to a possible replacement of $P_1$ to $P_2$ and $W_1$ to $W_2$, either
\begin{align}  \la{rp11}  &P_1=\nu \circ z^sR^n(z) \circ \sigma_1^{-1}, & &W_1=
\sigma_1 \circ z^n \circ \mu \\ 
&\la{rp11+} P_2=\nu \circ z^n \circ \sigma_2^{-1},& 
&W_2=\sigma_2 \circ
z^sR(z^n) \circ \mu, \end{align} 
where $R$ is a polynomial and $s\geq 0$, or
\begin{align} \la{rp21} &P_1=\nu \circ T_m \circ \sigma_1^{-1}, &
&W_1=\sigma_1 \circ T_n \circ \mu, \\
&\la{rp21+} P_2=\nu \circ T_n \circ \sigma_2^{-1}& &W_2=\sigma_2 \circ
T_m\circ \mu, \end{align} 
where $T_n, T_m$ are the Chebyshev polynomials. \ \qed
\et

Note that 
condition \eqref{ass0} implies that $\GCD(s,n)=1$ in \eqref{rp11}, \eqref{rp11+}.
In particular, 
the inequality $s>0$ holds whenever $n>1$.

For the reader convenience, in conclusion of this section we will make several comments concerning the proof of 
Theorem \ref{r2}.
First, if for given polynomials $P_1,$ $P_2$ the equality 
\be \la{erer}\GCD(\deg P_1, \deg P_2)=1\ee holds, then 
polynomials $Q_1,$ $Q_2$ satisfying equality \eqref{ura} exist if and only if the algebraic curve 
\be \la{cur} P_1(x)-P_2(y)=0\ee has genus zero, 
since condition \eqref{erer}
implies that \eqref{cur} is irreducible
and, in case if $g=0$, may be parametrized by polynomials.
Furthermore, if \eqref{cur} has genus zero, and $W_1,$ $W_2$ 
are $\widetilde W_1,$ $\widetilde W_2$ are two polynomial 
parametrizations of \eqref{cur} such that 
$$\deg W_1=\deg \widetilde W_1= \deg P_2, \ \ 
\deg W_2=\deg \widetilde W_2=\deg P_1,$$ then 
\be \la{ebi} P_1\circ W_1\sim P_1\circ  \widetilde W_1, \ \ \ \ P_2\circ W_2\sim P_2\circ \widetilde W_2\ee (for proofs of
the above statements see e.g. \cite{plau}, Section 2-4). 

Finally, note that the genus of \eqref{cur} depends on branch data of $P_1$ and $P_2$ only 
(see e.g. \cite{f1} or \cite{plau}, Lemma 8.2) and a direct though laborious analysis of the corresponding formula implies that the only possible passports of $P_1$ 
and $P_2$ for which $g=0$ are as in Theorem \ref{r2}. On the other hand, it is clear that 
$W_1,$ $W_2$ given in Theorem \ref{r2} provide parametrizations of the corresponding curves.


\section{Polynomial solutions of 
$P_1\circ W_1=P_2\circ W_2=\dots =P_r\circ W_r$
}
\subsection{Reduction to the case of coprime degrees. }
Similarly to the description of solutions of equation \eqref{ura} the description of solutions  
of the equation 
\be \la{ura2} P_1\circ W_1=P_2\circ W_2=\dots =P_r\circ W_r, \ee
where $P_i,W_i,$ $1\leq i \leq r,$ are polynomials 
of degrees $p_i,$ $w_i,$ $1\leq i \leq r,$ respectively, reduces to the case 
where \be \la{cop} \GCD(p_1,p_2,\dots ,p_r)=1,  \ \ \ \GCD(w_1,w_2,\dots ,w_r)=1.\ee

\bt \la{nr1} Let $P_i,W_i,$ $1\leq i \leq r,$ be polynomials such that \eqref{ura2} holds.
Then there exist polynomials
$U,$ $V,$ and $\widetilde P_i,\widetilde W_i,$ $1\leq i \leq r,$ where
$$\deg U=\GCD(p_1,p_2,\dots ,p_r),  \ \ \ \deg V=\GCD(w_1,w_2,\dots ,w_r),$$
such that
$$P_i=U\circ \widetilde P_i, \ \  \ W_i=\widetilde W_i\circ V, \ \ \ 1\leq i \leq r,$$
and $$\widetilde P_1\circ \widetilde W_1=\widetilde P_2\circ \widetilde W_2=\dots =\widetilde P_r\circ \widetilde W_r.\ \ \ \Box$$ 
\et
\pr The proof is essentially the same as in the case where $r=2$ since if
$B_i\in I_{d_i},$ $1\leq i \leq r,$ are blocks containing an element $x\in X$, then 
$\cap_{i=1}^r B_i$ is a block which coincides
with a residue class modulo $\LCM(d_1,d_2, \dots d_r)$, and  
$$\cap_{i=1}^r V_{d_i}=V_{\GCD(d_1,d_2, \dots d_r)}.\ \ \ \Box$$

\subsection{Proof of Theorem \ref{os1}}
The proof is by induction on $r$. For $r=2$ the statement follows from Theorem \ref{r2}. Assume now that the statement is true for $r-1$ and show that then it is true for $r\geq 3.$    
Show first that at least one $P_i,$ $1\leq i \leq r,$ 
is linearly equivalent either to a Chebyshev polynomial
or to a power. 

For $i,$ $1\leq i \leq r,$ set
$$x_i=\GCD(p_1,p_2,\dots p_{i-1},p_{i+1}, \dots p_r).$$
If at least one $x_i,$ $1\leq i \leq r,$ is equal to one, then 
the equality 
\be \l{per} P_1\circ W_1=P_2\circ W_2=\dots =P_{i-1}\circ W_{i-1}=P_{i+1}\circ W_{i+1}=\dots =P_{r}\circ W_{r}\ee
by the induction assumption implies that at least one $P_j,$ $1\leq j \leq r,$
$j\neq i,$  
is linearly equivalent either to a Chebyshev polynomial
or to a power. So, we may assume that \be \la{gep} x_i>1, \ \ \ 1\leq i\leq r.\ee 
Observe that condition \eqref{cop} implies that at least one of numbers $p_i,$ $1\leq i \leq r,$ 
is odd and without loss of generality we may assume that this is $p_r$. This 
implies that the numbers $x_i,$ $1\leq i\leq r-1,$ also are odd.

By Theorem \ref{nr1} there exist a polynomial $X_r,$ $\deg X_r=x_r,$ and polynomials $\widetilde P_i,$ $1\leq i \leq r-1,$ such that 
$$P_i=X_r\circ \widetilde P_i, \ \ \ 1\leq i \leq r-1,$$ 
and \be \la{ebbs} \widetilde P_1\circ W_1=\widetilde P_2\circ W_2=\dots =
\widetilde P_{r-1}\circ W_{r-1}.\ee 
Moreover, by the induction assumption at least one of polynomials $\widetilde P_i,$ $1\leq i \leq r-1,$ 
is linearly equivalent either to a Chebyshev polynomial
or to a power, and without loss of generality we may assume that 
this is $\widetilde P_1.$

Since \eqref{cop} implies that $\GCD(x_r,p_r)=1$ it follows from Theorem \ref{r1}
and Theorem \ref{r2} applied to the equality
$$P_r\circ W_r=X_r\circ (\widetilde P_1\circ W_1)$$  
that either $$P_r\sim T_{p_r},\ \ \ X_r\sim T_{x_r},$$ or $$P_r\sim z^{p_r},\ \ \ X_r\sim z^{s}R^{p_r}(z),$$ 
or 
$$P_r\sim  z^{s}R^{x_r}(z),\ \ \ X_r \sim z^{x_r},$$ 
where $R$ is a polynomial and $s\geq 0.$
Clearly, in the first two cases $P_r$  
is linearly equivalent either to a Chebyshev polynomial
or to a power and hence the statement is true. Therefore, we may assume that $X_r\sim z^{x_r}.$

In the similar way we may find polynomials $X_{r-1},$ $\deg X_{r-1}=x_{r-1},$ and $\widehat P_i,$ $1\leq i \leq r,$ $i\neq r-1,$ such that 
$$P_i=X_{r-1}\circ \widehat P_i, \ \  \ 1\leq i \leq r, \ \ \ i\neq r-1,$$  
and \be \la{ebbs2} \widehat P_1\circ W_1=\widehat P_2\circ W_2=\dots =
\widehat P_{r-2}\circ W_{r-2}=\widehat P_{r}\circ W_{r}.\ee
Furthermore, applying Theorem \ref{r1}
and Theorem \ref{r2} to the equality
$$P_{r-1}\circ W_{r-1}=X_{r-1}\circ (\widehat P_1\circ W_1)$$ we conclude as above that $X_{r-1}\sim z^{x_{r-1}}$, unless 
$P_{r-1}$  
is linearly equivalent either to a Chebyshev polynomial
or to a power.

Consider now the equality \be \la{xry} P_1= X_{r}\circ \widetilde P_1 =X_{r-1}\circ \widehat P_1\ee and show that if
\be \la{simm} X_r\sim z^{x_r}, \ \ \ X_{r-1}\sim z^{x_{r-1}},\ee and 
$\widetilde P_1$ is linearly equivalent either to a Chebyshev polynomial
or to a power, then $P_1$ is linearly equivalent to a power. 

First, observe that condition \eqref{cop} implies the equality \be \la{gep1} \GCD(x_r,x_{r-1})=1.\ee
In particular, at least one of the numbers $x_r,$ $x_{r-1}$ is greater than two. 
Therefore, since $z^n$ is not linearly equivalent to $T_n$ unless $n=2,$ 
it follows from Theorem \ref{r1}
and Theorem \ref{r2} 
applied to equality \eqref{xry} that there exist polynomials $W$, $\widetilde R$ and polynomials
$\alpha,$ $\beta,$ $\gamma$ of degree one such that 
either
\begin{align} &X_r=\alpha\circ z^{x_r}\circ \beta,& &\widetilde  P_1=\beta^{-1}\circ z^{s}\widetilde R(z^{x_r})\circ W,\\ \la{ega1}
&X_{r-1}=\alpha\circ z^{s}\widetilde R^{x_r}(z)\circ \gamma,&  &\widehat P_1=\gamma^{-1}\circ z^{x_r}\circ W,\end{align} 
or 
\begin{align}\la{ega2} &X_r=\alpha\circ z^{s}\widetilde R^{x_{r-1}}(z)\circ \beta,& &\widetilde  P_1=\beta^{-1}\circ z^{x_{r-1}}\circ W,\\ 
&X_{r-1}=\alpha\circ z^{x_{r-1}}\circ \gamma,&  &\widehat P_1=\gamma^{-1}\circ z^{s}\widetilde R(z^{x_{r-1}})\circ W.\end{align} 
Note that \eqref{gep} and \eqref{gep1} imply the inequality $s> 0.$

Observe now that if a polynomial $P$ of the form $z^s R^{n}(z),$ where $n>1$ and $s>0$, is linearly equivalent to a power, then $R$ is a monomial. Indeed, since a power has a unique  
critical point, the inequality $n>0$ implies that $R$ has at most one zero. Furthermore, since the multiplicity of the unique critical point of a power coincides with its degree  
it follows from $s>0$ that the unique zero of $R$ coincides with the origin. Therefore, it follows from \eqref{simm}, \eqref{ega1}, \eqref{ega2} that without loss of generality we may assume that 
\begin{align}\la{es1} &X_r=\alpha\circ z^{x_r}\circ \beta,& &\widetilde P_1=\beta^{-1}\circ z^{x_{r-1}}\circ W,\\  &X_{r-1}=\alpha\circ z^{x_{r-1}}\circ \gamma,& &\widehat P_1=\gamma^{-1}\circ z^{x_{r}}\circ W, 
\end{align}
where $W$ is a polynomial and 
$\alpha,$ $\beta,$ $\gamma$ are polynomials of degree one.

If $\widetilde P_1$ is linearly equivalent to a power, then it follows from the second equality in \eqref{es1} by the chain rule that the only critical value of $W$ is 
zero implying that
$W=z^t\circ \omega$ for some polynomial of degree one $\omega$ and $t\geq 0$. Therefore, in this case 
$P_1=X_r\circ \widetilde P_1$ is linearly equivalent to a power. 
On the other hand, 
the above assumptions yield that $\widetilde P_1$ may not be  
linearly equivalent to a Chebyshev polynomial for otherwise 
Corollary \ref{c1} applied to the second equality in \eqref{es1}
would imply that $z^{x_{r-1}}$ is linearly equivalent to 
$T_{x_{r-1}}$ in contradiction with the 
assumption that 
$x_{r-1}$ is an odd number greater than one. 
 
In order to prove that at least one $W_i,$ $1\leq i \leq r,$ 
is linearly equivalent either to a Chebyshev polynomial
or to a power we use similar arguments. Namely, 
for $i,$ $1\leq i \leq r,$ define
$$ y_i=\GCD(w_1,w_2,\dots w_{i-1},w_{i+1}, \dots w_r).$$
As above, if at least one $y_i,$ $1\leq i \leq r,$ is equal to one, then equality \eqref{per} 
by the induction assumption implies that at least one $W_j,$ $1\leq j \leq r,$
$j\neq i,$  
is linearly equivalent either to a Chebyshev polynomial
or to a power. So, we may assume that $y_i>1$ for all $i,$ $1\leq i\leq r.$
Furthermore, we may assume that $w_r$ and $y_i,$ $1\leq i\leq r-1,$ are odd.

Using Theorem \ref{nr1} we conclude that there exist a polynomial $Y_r,$ $\deg Y_r=y_r,$ and polynomials $\widetilde W_i,$ $1\leq i \leq r-1,$ such that 
$$W_i=\widetilde W_i\circ Y_r, \ \ \ 1\leq i \leq r-1,$$ 
and \be  P_1\circ \widetilde W_1=P_2\circ \widetilde W_2=\dots =
P_{r-1}\circ\widetilde  W_{r-1},\ee 
where by the induction assumption we may assume that 
$\widetilde W_1$
is linearly equivalent either to a Chebyshev polynomial
or to a power. Furthermore, since \eqref{cop} implies that $\GCD(y_r,w_r)=1$ it follows from Theorem \ref{r1} and Theorem \ref{r2} applied to the equality
$$(P_1\circ  \widetilde W_1)\circ Y_r=P_r\circ W_r$$ 
that $W_r$ 
is linearly equivalent either to a Chebyshev polynomial
or to a power unless $Y_r\sim z^{y_r}$.

Continuing arguing as above we reduce the proof of the theorem to the analysis of the equality  
\be \la{xry1} W_1= \widetilde W_1 \circ Y_{r}=\widehat W_1\circ Y_{r-1}, \ee where  
\be \la{simm+} Y_r\sim z^{y_r}, \ \ \ Y_{r-1}\sim z^{y_{r-1}},\ee 
$\widetilde W_1$ is linearly equivalent to a Chebyshev polynomial
or to a power, and $\widehat W_1$ is a polynomial.

Observe now that if a polynomial of the form $z^s R(z^{n}),$ where $n>1,$ $s>0,$
is linearly equivalent to a power, then $R$ is a monomial. Indeed,
comparing the coefficients of $z^{n-1}$ of both parts of  
the equality $$z^s R(z^{n}) =\mu \circ z^l\circ \nu,$$  
we conclude that 
$\nu(0)=0.$ It follows now from $s>0$ that 
$\mu(0)=0$ implying that $R$ is a monomial. Therefore, applying 
Theorem \ref{r1} and Theorem \ref{r2} to equality \eqref{xry1} 
and arguing as in the analysis of equality \eqref{xry} 
we 
conclude that there exist a polynomial $W$ and polynomials
$\alpha,$ $\beta,$ $\gamma$ of degree one such that 
\begin{align} \la{es1+} &\widetilde  W_1=W\circ z^{y_{r-1}}\circ \beta,& &Y_r=\beta^{-1}\circ z^{y_{r}}\circ \alpha,\\  &\widehat W_1=W\circ z^{y_{r}}\circ \gamma,&  &Y_{r-1}=\gamma^{-1}\circ z^{y_{r-1}}\circ \alpha. 
\end{align}

If $\widetilde W_1$ is linearly equivalent to a power, then 
the first equality in \eqref{es1+} implies that 
$W$ has a unique critical value and that the corresponding critical point is zero
for otherwise $\widetilde W_1$ 
would have more than one critical point. Therefore, 
$W=\omega\circ z^t$ for some polynomial of degree one $\omega$ and $t\geq 0$ implying that 
$W_1=\widetilde W_1\circ Y_r$ is linearly equivalent to a power. 
On the other hand, $\widetilde W_1$ may not be  
linearly equivalent to a Chebyshev polynomial since otherwise Corollary \ref{c1} applied to the first equality in \eqref{es1+} would imply that 
$z^{y_{r-1}}\sim T_{y_{r-1}}$ in contradiction with the assumption that $y_{r-1}$ is an odd number greater than one.  \qed

\subsection{\la{appp} Double decompositions involving Chebyshev polynomials or powers}
By Theorem \ref{os2} if a polynomial $P$ has several ``coprime" compositional right factors, then
one of this factors is linearly equivalently either to a Chebyshev or to a power. In this subsection we describe the form 
of other right factors of such $P$ and the form of $P$ itself.

\bl \la{c2} Let $P,P_1,W_1,W_2$ be polynomials  
such that
$$P=P_1\circ z^n=P_2\circ W_2.$$ 
Then there exist 
polynomials $R,U$ and 
a polynomial 
$\sigma$ of degree one such that 
\be \la{go1} W_2=  \sigma\circ z^sR(z^n),\ \ \ P=U\circ z^{sn/e}R^{n/e}(z^n),\ee where $s\geq 0$ and $e=\GCD(n,\deg W_2).$\el

\pr 
Observe first that without loss of generality we may assume  
that \be \la{ass00} \GCD(\deg P_1,\deg P_2)=1, \ \  \GCD(n,\deg W_2)=1.\ee
Indeed, by Theorem \ref{r1}
that there exist polynomials $ A,$ $B,$ $C,$ $D,$ $U,$ $V$ where $\deg U=\GCD(\deg P_1,\deg P_2),$ $\deg V=e,$ 
such that $$ P_1=U \circ A, \ \  \ \ z^n=C\circ V, \ \ \ \ P_2=U \circ B,\ \ \ \ 
W_2=D\circ  V, $$ 
and \be \la{laf} A\circ C=  B\circ D.
\ee 
Furthermore, it follows from the first part of Corollary \ref{c1} that without loss of generality we may assume that 
$$C=z^{n/e}, \ \  \ V= z^{e}.$$ 
Denote the polynomial defined by equality \eqref{laf} by $\widetilde P$.
If the proposition is true under assumption \eqref{ass00}, then $$D=  \sigma\circ z^{l}R(z^{n/e}),\ \ \ \widetilde P=z^{ln/e}R^{n/e}(z^{n/e}),$$ where $\GCD(l,n/e)=1.$
Therefore, since $$P=U\circ \widetilde P\circ z^e, \ \ \ W_2=D\circ z^e,$$ equalities \eqref{go1} hold with $s=le.$ 

In order to prove Lemma \ref{c2} 
under assumption \eqref{ass00} 
one can use Theorem \ref{r2}. However, 
a shorter way is to observe 
that the equality $P=P_1\circ z^n$ implies that for any primitive $n$th root of unity $\v$ the equality $P(\v z)=P(z)$ holds. Therefore,
$$P_2\circ W_2=P_2\circ (W_2 \circ \v z)$$ and applying to this equality Lemma \ref{c0} we conclude that 
\be \la{los} W_2=\mu\circ W_2 \circ \v z\ee for some polynomial $\mu $ of degree one. It follows now from the comparison of coefficients of parts of \eqref{los} that $W_2=z^{s}R(z^n)+\alpha$ for some $R\in \C[z],$ $\alpha\in \C,$ and $s\geq 0.$ 

Further, observe that for given $W_1$, $W_2$ polynomials $P_1,$ $P_2$ such that equalities \eqref{ura} and \eqref{ass0} hold are defined in a unique 
way up to the change $P_1\rightarrow \mu \circ P_1,$ $P_2\rightarrow \mu \circ P_2,$
where $\mu$ is a polynomial of degree one.
Indeed, let $\widetilde P_1,$ $\widetilde P_2$ be another such polynomials. Then 
\eqref{ura} is also satisfied for 
$$\widehat  P_1=P_1-\lambda \widetilde P_1, \ \ \ \widehat P_2=P_2-\lambda \widetilde P_2,$$ where  $\lambda$ is any complex number,
and hence choosing appropriate $\lambda$ 
we may obtain a pair 
$\widehat P_1,$ $\widehat P_2$ of polynomials 
such that \be \la{iiuuo} \widehat P_1\circ W_1=\widehat P_2\circ W_2\ee
and \be \la{plm} \deg \widehat P_1<\deg P_1=\deg W_2.\ee On the other hand, it is easy to see 
comparing the  
leading terms of parts of equality \eqref{iiuuo} and taking into account the equality $\GCD(W_1,W_2)=1$
that 
\eqref{plm} may not be satisfied 
unless $\widehat P_1\equiv\widehat P_2\equiv c$ for some $c\in \C$. 

Therefore, since the equality 
\be \la{iqw} P_1\circ z^n=P_2\circ  (\sigma\circ z^{s}R(z^n))\ee is clearly satisfied for 
$$P_1=z^{s}R^n(z), \ \ \  P_2=z^n\circ \sigma^{-1},$$ we conclude that 
$P$ has the form indicated in \eqref{go1}. \qed

\bl \la{c3} Let $P,P_1,W_1,W_2$ be polynomials 
such that \be \la{sds} P=P_1\circ T_n=P_2\circ W_2\ee and 
$n\nmid \deg W_2$. Then  
there exist a polynomial $U$ and a polynomial $\sigma$ of degree one 
such that 
either
\be \la{pia11} W_2=  \sigma\circ T_m, \ \  P=U\circ T_{t}, \ee 
where $t=\LCM(n,m),$
or
\be \la{pia12} W_2= \sigma\circ zS(z^2) \circ T_{n/2}, \ \  
P=U\circ z^2S^2(z^2) \circ T_{n/2},\ee for some polynomial $S$. 
\el

\pr 
Using the second part of Corollary \ref{c1} it is easy to show in the same way as in the proof 
of Lemma \ref{c2} that without loss of generality we may assume that condition \eqref{ass00} holds. Furthermore, the condition $n\nmid \deg W_2$ implies that $n\geq 2.$ 

If $n=2$, then, since $T_2=\theta \circ z^2$, where $\theta=2z-1,$ 
the lemma follows from Lemma \ref{c2}  
taking into account that $s=1$ in formulas \eqref{go1} in view of the condition $n\nmid \deg W_2$.

Assume now that $n> 2$ and apply Theorem \ref{r2} to equality \eqref{sds}. If equalities \eqref{rp21}, \eqref{rp21+} hold, then the statement obviously is true. Otherwise, taking into account that $z^n$ and $T_n$ are not linearly equivalent for $n>2,$ we conclude that there exist polynomials $\sigma_1,$ $\sigma_2,$ $\nu,$ $\mu$ of degree one  
such that
\begin{align} \la{kor0} &P_1=\nu \circ z^{n_1} \circ\sigma_1^{-1},& &T_n=\sigma_1 \circ z^{s_1}R_1(z^{n_1}) \circ \mu, \\ \la{kor}
&P_2= \nu \circ
z^{s_1}R_1^{n_1}(z) \circ \sigma_2^{-1}, &
&W_2= \sigma_2 \circ z^{n_1} \circ \mu,\end{align} 
where $R_1$ is a non-constant polynomial and $\GCD(s_1,n_1)=1$. 

If $n_1=1,$ then the lemma is clearly true. So, assume that $n_1>1.$ It is not hard to see that 
if $\zeta$ is a critical point of $z^{s_1}R_1(z^{n_1})$, then for any $i,$ $1\leq i \leq n_1,$ the number
$\v^i\zeta,$ where $\v$ is an $n_1$th primitive root of unity, also is a critical point.
On the other hand, since 
all critical points of $T_n$ are
on the real line, the equality
\be \la{kon}
T_n=\sigma_1 \circ z^{s_1}R_1(z^{n_1}) \circ \mu
\ee
implies that all critical points of 
$$z^{s_1}R_1(z^{n_1})=\sigma_1^{-1} \circ T_n\circ \mu^{-1}$$ are on the line
$\mu\{\mathbb R\}$. This implies easily that if \eqref{kon} holds, then $n_1=2$ and $\mu =\pm z.$ Therefore, $W_2=\sigma \circ T_2,$ where $\sigma=\sigma_2 \circ \,\theta^{-1}$ and we can finish the proof as in Lemma \ref{c2} observing that $T_n$ and $\sigma \circ T_2$ parametrize the curve $$T_2(x)-(T_n\circ \sigma^{-1})(y)=0.\ \ \Box$$

\vskip 0.2cm
\noindent {\bf Remarks.} Notice that the proofs of Lemma \ref{c2} and Lemma \ref{c3} given above actually describe not only possible forms of $W_2$ but also possible forms of $P_1$ and $P_2$. Notice also that one can obtain similar descriptions of solutions of \eqref{ura} in the case where a {\it left} compositional factor of $P$ 
is a Chebyshev polynomial or to a power. For this purpose one can use Theorem \ref{r2} or the genus formula for curve \eqref{cur}. 

\pagebreak
Finally, notice that for small $r$ one can try to obtain a more detailed description of solutions of \eqref{ura2} in the spirit of the second Ritt theorem.
For examples, one can show that 
any solution of the equation 
$$ z^n\circ A=B \circ z^m=U\circ V $$
has the form 
$$U=z^{\frac{r_2m}{d_2}}R_2^n(z^{\frac{m}{d_2}}), \ \ \ 
V=z^{\frac{r_1n}{d_1}}R_1^{\frac{n}{d_1}}(z^m),
$$
$$A=z^{\frac{r_1r_2m}{d_1d_2}}R_1^{\frac{r_2m}{d_1d_2}}(z^m)
R_2(z^{\frac{r_1nm}{d_1d_2}}R_1^{\frac{mn}{d_1d_2}}(z^m))), \ \ \ 
B=z^{\frac{r_1r_2n}{d_1d_2}}R_1^{\frac{n}{d_1}}
R_2^n(z^{\frac{r_1n}{d_1d_2}}R_1^{\frac{mn}{d_1d_2}}),
$$
where $R_1,$ $R_2$ are polynomials, $\GCD(r_1,m)=1,$ $\GCD(r_2,n)=1,$ and
$d_1d_2=\GCD(n,m).$ 
However, similar results seem not to have valuable applications.

\section{\la{pp2a} Proof of Theorem \ref{os2}}
Recall that owing to the result of \cite{pm} in order to prove Theorem \ref{os2} we only must show that if polynomials $P$ and $Q=\sum_{i=1}^r Q_i$ satisfy 
\be \la{pqp}
P=P_i\circ W_i, \ \ \
Q_i=V_i\circ W_i, \ \ \  W_i(a)=W_i(b), \ \ \ 1\leq i \leq r,
\ee
for some polynomials $P_i, V_i, W_i,$ $1\leq i \leq r$, and $a\neq b$, then one of conclusions 1)-4) of Theorem \ref{os2} holds. 
Note that for any polynomials $\nu$, $\mu_i$, $1\leq i \leq s,$ of degree one equalities \eqref{pqp} imply 
similar equalities 
for polynomials 
$$\widetilde P_i=P_i\circ \mu_i^{-1}, \ \ \ \widetilde V_i=V_i\circ \mu_i^{-1},\ \ \ \widetilde W_i=\mu_i \circ W_i\circ \nu, \ \ \  1\leq i \leq s,$$
and vice versa.
We often will 
use this property in order to simplify calculations.

The proof splits into two parts. First,  
we will prove the following proposition. 

\bp \la{ese} Let $P$, $P_i,$ $W_i,$ $1\leq i \leq r$, be polynomials and $a, b$ be complex numbers such that  
the  
equalities 
\be \la{pqpq}
P=P_i\circ W_i, \ \ \ W_i(a)=W_i(b), \ \ \ 1\leq i \leq r,
\ee hold for some $r>3$. Then there exists a pair of distinct indices $i_1,i_2,$ $1\leq i_1,i_2 \leq r,$ 
such that \be \la{iiuu} W_{i_1}=\widetilde W_{i_1}\circ Z, \ \ \  W_{i_2}=\widetilde W_{i_2}\circ Z, \ \ \ Z(a)=Z(b)\ee for some polynomials $\widetilde W_{i_1},$ $\widetilde W_{i_2},$ and $Z.$
\ep

Proposition \ref{ese} reduces the proof of the theorem to the case where $r\leq 3$. Indeed, if $r>3,$ then
the number of reducible solutions in the representation $Q=\sum_{i=1}^r Q_i$ always may be reduced by one since 
the sum of two reducible
solutions $$Q_{i_1}=V_{i_1}\circ W_{i_1}, \ \ \  Q_{i_2}=V_{i_2}\circ W_{i_2},$$ may be replaced by 
the unique reducible 
solution $$(V_{i_1}\circ \widetilde W_{i_1}+V_{i_2}\circ \widetilde W_{i_2})\circ  Z.$$  
The second part of the proof consists in the analysis of condition \eqref{pqpq} in the case where $r\leq 3.$
 
We start by proving the following technical lemma.

\bl \la{skun} Let $T_{m_1},$ $T_{m_2},$ $T_{m_3}$ be the Chebyshev polynomials and $a$, $b$ be complex numbers.
\vskip 0.1cm
\noindent a) Assume that 
\be \la{uuiioo} T_{m_1}(a)=T_{m_1}(b), \ \ \ \ \ \ T_{m_2}(a)=T_{m_2}(b), \ \ \ \ \ \ T_{m_3}(a)=T_{m_3}(b).\ee
Then there exists a pair of distinct indices $i_1,i_2,$ $1\leq i_1,i_2 \leq 3,$ 
such that for $l=\GCD(m_{i_1},m_{i_2})$ the equality $T_l(a)=T_l(b)$ holds.
\vskip 0.1cm
\noindent b) Assume that 
\be \la{volk} T_{m_1}(a)=0, \ \ \ T_{m_2}(a)=0,\ee 
where $m_1,$ $m_2$ are odd numbers such that $\GCD(m_1,m_2)=1$. Then 
$a=0.$
\vskip 0.1cm
\noindent c) Assume that 

\be \la{volk2} T_{m_1}(a)=-T_{m_1}(b), \ \ \ T_{m_2}(a)=-T_{m_2}(b),\ee
where $m_1,$ $m_2$ are odd numbers such that $\GCD(m_1,m_2)=1$. Then 
either $a=-b,$ or 
$T_{m_1m_2}(a)=\pm 1.$ 
\el

\pr Choose $\alpha, \beta \in \C$ such that $\cos \alpha= a,$
$\cos \beta =b.$ Then equalities \eqref{uuiioo}
imply the equalities
\be \la{10} m_1 \alpha = \pm m_1 \beta + 2\pi k_1, \ \ \
m_2 \alpha = \pm m_2 \beta+ 2\pi k_2, \ \ \ 
m_3 \alpha = \pm  m_3 \beta+ 2\pi k_3, \ee where $k_1,k_2,k_3\in \Z.$
Clearly, the signs in at least two equalities
\eqref{10} coincide. To be definite suppose that they coincide 
in the first two equalities and choose $u,v\in \Z$ such that $um_1+vm_1=l,$ where $l={\rm GCD}(m_1,m_2).$ 
Multiplying now the 
first equality in \eqref{10} by $u$ and adding 
the second equality multiplied by $v$ we 
see that $l\alpha= \pm l\beta + 2\pi k_4,$ where $k_4\in \Z.$ Therefore,
$\cos \, l\alpha=\cos \, l\beta$ implying that 
$T_{l}(a)=T_{l}(b)$. 

Further, equalities \eqref{volk} imply the equalities   
\be \la{be} m_1\alpha =\pi/2 +\pi k_1, \ \ \ m_2\alpha =\pi/2 +\pi k_2, \ \ \ k_1,k_2\in \Z.\ee 
Multiplying the 
first equality in \eqref{be} by $u$ and adding 
the second equality multiplied by $v$, 
where $u,v$ satisfy 
\be \la{pes} um_1+vm_2=1,\ee  we 
see that \be \la{pes2} \alpha=(u+v)\pi/2+ \pi k_3, \ \ \ k_3\in \Z.\ee 
Moreover, since $m_1$, $m_2$ are odd, equality \eqref{pes} implies that the numbers $u,$ $v$ have different parity. Therefore, \eqref{pes2} implies that 
$a=\cos \alpha=0$. 

Finally, equalities \eqref{volk2} imply the equalities  
\be \la{cos1} m_1\alpha=\pi \pm m_1\beta +2\pi k_1, \ \  \ m_2\alpha=\pi \pm m_2\beta +2\pi k_2, \ \  \ k_1, k_2\in \Z.\ee 
If the signs in equalities \eqref{pes} are the same, then
$$\alpha=(u+v)\pi\pm \beta+2\pi k_3, \ \ \ k_3\in \Z,$$ 
where $u,v$ satisfy \eqref{pes}. Since the numbers $u,$ $v$ have different parity, this 
implies that $a=-b.$ On the other hand, if the signs in \eqref{cos1} are opposite, then  
multiplying the first equality in \eqref{pes} by $m_2$ and adding to the second equality multiplied by $m_1$ we conclude that
$$2m_1m_2\alpha=\pi (m_1+m_2) +2\pi k_3, \ \  \ k_3\in \Z.$$ Owing to the oddness of $m_1$ and $m_2$ this implies that 
$T_{2m_1m_2}(a)=1.$ Therefore, since $T_{2m_1m_2}=T_2\circ T_{m_1m_2}$ and $T_2=2z^2-1$, the equality  
$T_{m_1m_2}(a)=\pm 1$ holds. \qed

\vskip 0.2cm
\noindent{\it Proof of Proposition \ref{ese}.} 
First of all observe that we may assume that 
\be \la{cop1}  \GCD(\deg W_1,\deg W_2,\dots ,\deg W_r)=1.\ee Indeed, 
if $$\GCD(\deg W_1,\deg W_2,\dots ,\deg W_r)=w>1,$$ then it follows from \eqref{pqpq} taking into account Theorem \ref{nr1} that there exist polynomials $\widehat P,$
$\widehat W_i,$  
$1\leq i \leq r,$ and $W,$ 
$\deg W=w,$ such that  
\be \la{st1} W_i=\widehat W_i\circ W, \ \ \ 1\leq i \leq r,\ee
and 
\be \la{st2}
\widehat P=P_i\circ \widehat W_i, \ \ \
\widehat W_i(W(a))=\widehat W_i(W(b)), \ \ \ 1\leq i \leq r.
\ee
Therefore, if 
the proposition is true under assumption \eqref{cop1} and $\widetilde W_{i_1},$ $\widetilde W_{i_2},$ $U$ are polynomials such that 
$$\widehat  W_{i_1}=\widetilde W_{i_1}\circ U, \ \ \  \widehat W_{i_2}=\widetilde W_{i_2}\circ U, \ \ \ U(W( a))=U(W(b)),$$ then equalities \eqref{iiuu} hold for the same $i_1,$ $i_2$ and 
$Z=U\circ W.$

By Theorem \ref{os1} 
equality \eqref{cop1} implies that 
we may assume that $W_1$ is equivalent either to $z^n$ or $T_n$.  
Furthermore, owing to the remark made in the beginning of this section,  
without loss of generality we may assume that $W_1$ is equal either to $z^n$ or 
$T_n$. Note that the equality $W_1(a)=W_1(b)$ implies that $n>1$.

\vskip 0.2cm
\noindent{\it Case 
$W_1=z^n$}. Since $P=P_1\circ z^n$ it follows from Lemma
\ref{c2} applied to the equality
\be \la{begi} P=P_1\circ z^n=P_2\circ W_2\ee
that $W_2=\sigma\circ z^{s}R(z^n)$, where $R$ is a polynomial and $\sigma$ is a polynomial of degree one.
In particular, this implies that if $\v$ is a primitive $n$th root of unity, then equality 
\be \la{los1} W_2=\mu\circ W_2 \circ \v z\ee
holds for some 
polynomial $\mu$ of degree one.

Further, it follows from 
Theorem \ref{r1} and Theorem \ref{r2} applied to the equality $$P_2\circ W_2=P_3\circ W_3$$ that without loss of generality we may assume that either 
\be \la{il} W_2=z^m\circ W, \ \ \ \ W_3=z^{s_1}R_1(z^m)\circ W,\ee where $R_1,$ $W$ are polynomials
and $\GCD(m,s_1)=1,$  
or
\be \la{il1} W_2=T_{m_1}\circ W, \ \ \ \ W_3=T_{m_2}\circ W,\ee where $W$ is a polynomial and $\GCD(m_1,m_2)=1.$
If $W(a)=W(b)$, then \eqref{iiuu} obviously holds. 
So, below we may assume that 
\be \la{nee} W(a)\neq W(b). \ee This implies in particular that $m>1$ in \eqref{il}.

Assume first that \eqref{il} has place. Then \eqref{los1} implies that 
$$z^m\circ W=(\mu\circ z^m)\circ (W\circ \v z)$$ and applying now Lemma \ref{c0}
to the last equality we conclude that there exists a polynomial $\nu$ of degree one
such that \be \la{ser} z^m=\mu\circ z^m\circ \nu^{-1},\ee 
and \be \la{ser1}W =\nu\circ W\circ \v z.\ee 
Since $m>1$, equality \eqref{ser} implies that $\nu(0)=0$
and the comparison of coefficients of the parts of \eqref{ser1} yields that \be \la{rf1} W=z^{s_2}R_2(z^n)\ee for some $s_2\geq 0$ and a polynomial $R_2$. Therefore, \be \la{rf2} W_2=z^{ms_2}R_2^m(z^n)\ee by \eqref{il}. 

Clearly, equalities $W_1(a)=W_1(b),$ $W_2(a)=W_2(b)$ imply that either 
the number $a^n=b^n$ is a root of $R_2$ or the equality $a^{ms_2}=b^{ms_2}$ holds.
In the first case equality \eqref{rf1} implies that $W(a)=W(b)$ while in the second one 
we conclude that $ a^t=b^t,$ where $t=\GCD(ms_2,n)$,  
implying that \eqref{iiuu} holds for $i_1=1,$ $i_2=2,$ and $Z=z^t.$

Assume now that \eqref{il1} holds. Then \eqref{los1} implies that
$$T_{m_1}\circ W=(\mu\circ T_{m_1})\circ (W\circ \v z)$$ and applying to this equality Lemma \ref{c0}
we conclude that there exists a polynomial $\nu$ of degree one
such that the equalities \be \la{ser+} T_{m_1}=\mu\circ T_{m_1}\circ \nu^{-1},\ee 
and \eqref{ser1} hold. It follows from equality 
\eqref{ser+} that $\nu$ transforms the set of critical points of $T_{m_1}$ to itself. Since all critical points of $T_{m_1}$ are real this yields easily that $\nu=\pm z.$ 
Furthermore, if $\nu=z$, then \eqref{ser1} implies that $W=R_2(z^n)$ for some
polynomial $R_2$ and it follows from $a^n=b^n$ that $W(a)=W(b).$
On the other hand, if 
$\nu=-z$, then it follows from \eqref{ser1} that
$W=z^{n/2}R_2(z^n)$, 
and $a^n=b^n$ implies that either $W(a)=W(b)$ or
\be \la{nee2} W(a)=-W(b).\ee

Since $W_2=T_{m_1}\circ W$, the equalities $W_2(a)=W_2(b)$ and \eqref{nee2}, taking into account equality \eqref{ioi}, imply  that either 
$m_1$ is even or \be \la{vol1} T_{m_1}(W(a))=T_{m_1}(W(b))=0.\ee Similarly, $W_3(a)=W_3(b)$ and \eqref{nee2} imply 
that either $m_2$ is even or \be \la{vol2} T_{m_2}(W(a))=T_{m_2}(W(b))=0.\ee 
If at least one of numbers $m_1,$ $m_2$, say $m_1$, is even, then by \eqref{ioi} 
there exists a polynomial $F$ such that
$T_{m_1}=F\circ z^2$ and hence $$W_2=T_{m_1}\circ z^{n/2}R_2(z^n)=F\circ z^nR_2^2(z^n)=F\circ zR_2^2\circ z^n$$ implying that \eqref{iiuu} holds for 
$i_1=1,$ $i_2=2,$ and $Z=z^n.$ On the other hand, if both $m_1,$ $m_2$ are odd, then it follows from 
\eqref{vol1}, \eqref{vol2} by Lemma \ref{skun}, b) that 
$W(a)=W(b)=0.$ 

Note that in the above proof we actually did not use the assumption $r>3$ but only the weaker assumption 
$r>2$.

\vskip 0.2cm
\noindent{\it Case $W_1=T_n.$} 
Observe first that if $n$ is a divisor of $\deg W_{j},$ $2\leq j\leq r,$ then the proposition is true since 
Corollary \ref{c0} applied to the equality
\be \la{net} P_1\circ T_n=P_j\circ W_j\ee
implies that $W_{j}=R\circ T_n$ for some polynomial $R$, and hence \eqref{iiuu} holds for $i_1=1,$ $i_2=j,$ and $Z=T_n.$ Otherwise, 
Lemma \ref{c3} applied to \eqref{net} implies that without loss of generality we may assume that 
$W_j,$ $2\leq	j \leq r,$ either is a Chebyshev polynomial or has the form $zR(z^2)\circ T_{n/2}$ for some polynomial $R.$ Furthermore, 
since $r>3$, at least two polynomials from the set $W_j,$ $2\leq	j \leq r,$ either are both Chebyshev polynomials or 
both have the form $zR(z^2)\circ T_{n/2}$. Therefore, 
without loss of generality we may assume that either 
\be \la{ocel} W_2=T_{m_1}, \ \ \ W_3=T_{m_2},\ee  
or 
\be \la{ocel2} W_2=zR_1(z^2)\circ T_{n/2}, \ \ \ W_3=zR_2(z^2)\circ T_{n/2}\ee
for some polynomials $R_1,$ $R_2.$ 

If \eqref{ocel} holds, then the proposition is true by Lemma \ref{skun}, a). 
So assume that \eqref{ocel2} has place.
In this case equalities \eqref{pqpq} 
imply the equalities 
\be \la{ii1}P_1\circ \widehat W_1=P_2\circ \widehat W_2=P_3\circ \widehat W_3,\ee 
and  \be \la{ii} \widehat W_1(\widehat a)=\widehat W_1(\widehat b), \ \ \  \widehat W_2(\widehat a)=\widehat W_2(\widehat b),
\ \ \  \widehat W_3(\widehat a)=\widehat W_3(\widehat b),\ee
where 
$$ \widehat W_1=T_2, \ \ \ \widehat W_2= zR_1(z^2), \ \ \ \widehat W_3= zR_2(z^2),$$
and 
$$\widehat a=T_{n/2}(a),\  \ \  \widehat b=T_{n/2}(b).$$
Since $W_1=T_2\sim z^2$, it is already proved that \eqref{ii1} and \eqref{ii} 
imply 
that there exists a pair of indices $i_1,i_2,$ $1\leq i_1,i_2 \leq r,$ 
such that $$\widehat  W_{i_1}=\widetilde W_{i_1}\circ U, \ \ \  \widehat W_{i_2}=\widetilde W_{i_2}\circ U, \ \ \ U(\widehat a)=U(\widehat b)$$ for some polynomials $\widetilde W_{i_1},$ $\widetilde W_{i_2},$ and $U$. Therefore, \eqref{iiuu} holds for the same $i_1,$ $i_2$ and 
$Z=U\circ T_{n/2}.$
\vskip 0.2cm

\noindent{\it Proof of Theorem \ref{os2}.} Recall that owing to Proposition \ref{ese}  
in order to prove Theorem \ref{os2} we only must show that if \eqref{pqp} holds for $r\leq 3$, then one of the cases 1)-4) listed in the formulation 
has place. If 
$r=1$, then clearly 1) holds. So, assume that $r>1.$ 
An argument similar to the one given in the beginning of Proposition
\ref{ese} shows that without loss of generality  
we may assume that equality \eqref{cop1} holds and either $W_1=z^n$ or 
$W_1=T_n$. Furthermore, as it was remarked above if $W_1=z^n$, then we may suppose that $r=2$, and applying Lemma \ref{c2}   
to equality \eqref{begi} we see that without loss of generality we may assume that 
$$P=U\circ z^{sn}R^n(z^n),\ \ \ W_2= z^sR(z^n),$$
where $U$, $R$ are polynomials, $n>1$, $s>0$ and $\GCD(s,n)=1.$
Moreover, since $\GCD(s,n)=1$, 
the equalities $W_1(a)=W_1(b),$ $W_2(a)=W_2(b)$
imply that the number $a^n=b^n$ is a root of $R$. Therefore, 
if $W_1=z^n,$ then the second case listed in Theorem \ref{os2} has place.

In the case where $W_1=T_n$ the number $r$ may be equal to 2 or to 3.
Further, 
the analysis given in the proof of Proposition \ref{ese} implies that 
in the first case without loss of generality we may assume that
either \be \la{fc} W_1=T_n, \ \ \ W_2=T_m,\ee or 
\be \la{sc} W_1=T_n, \ \ \ W_2=zR(z^2)\circ T_{n/2},\ee 
while
in the second case   
we may assume that
\be \la{tc} W_1=T_n, \ \ \ W_2=T_m, \ \ \ W_3=zR_1(z^2)\circ T_{n/2},\ee where $R_1$ is a polynomial.
Furthermore, in the last case without loss of generality we may assume that 
\be \la{psol} W_3\neq \sigma\circ T_l\ee for a Chebyshev polynomial $T_l$ and a polynomial $\sigma$ of degree one.

If \eqref{fc} holds, then it follows from Lemma \ref{c3} applied to the equality 
\be \la{pezd} P=P_1\circ T_n=P_2\circ T_m\ee that the third 
case listed in Theorem \ref{os2} has place. On the other hand, if \eqref{sc} holds, then 
since $n/2$ is a divisor of both $\deg W_1$ and $\deg W_2$ it follows from \eqref{cop1} that $n/2=1$. Therefore, $W_1\sim z^2$ and 
hence by the above argument the second case listed in Theorem \ref{os2} has place. 

Consider finally 
the case where \eqref{tc} holds. 
Observe first that owing to  \eqref{psol}
Lemma \ref{c3} applied to the equality $$P_2\circ T_m=P_3\circ W_3$$ implies
that
\be \la{bu} W_3=z R_2(z^2)\circ T_{m/2}\ee 
for some polynomial $R_2$. In particular, this implies that $m$ is even.
Further, since $n/2$ divides both $\deg W_1$ and $\deg W_3$ it follows from \eqref{cop1} that $\GCD(n/2,m)=1$. 
Similarly, \eqref{bu} implies $\GCD(n,m/2)=1.$ This yields that $n/2$ and $m/2$ are 
odd and $\GCD(n/2,m/2)=1.$ 

Since $\GCD(n,m)=2,$ Lemma \ref{c3} applied to equality \eqref{pezd}
implies that 
\be \la{we2} P=V\circ T_{nm/2},\ee
where $V$ is a polynomial. 
Applying now Lemma \ref{c3} to the equation $$P=V\circ T_{nm/2}=P_3\circ W_3$$
and taking into account \eqref{psol} we conclude that \be \la{lak} W_3=z R(z^2)\circ T_{nm/4}, \ \ \ P=U\circ z^2 R^2(z^2)\circ T_{nm/4}\ee 
for some polynomials $R$ and $U$.

If $T_{n/2}(a)=T_{n/2}(b)$,
then we may replace $Q_1+Q_3$ by 
$$(V_{1}\circ T_2+V_{3}\circ zR_1(z^2))\circ  T_{n/2}$$ 
and hence $Q$ is a sum of only two reducible solutions. So, we may assume that
$T_{n/2}(a)\neq T_{n/2}(b)$ implying by 
$T_{n}(a)= T_{n}(b)$ that $T_{n/2}(a)=-T_{n/2}(b).$  
A similar argument shows that we may assume that $T_{m/2}(a)=-T_{m/2}(b)$.
Finally, if $a=-b,$ then we may replace $Q_1+Q_2$ by 
$$(V_{1}\circ T_n+V_{2}\circ T_m)\circ  T_{2},$$ so assume that $a\neq -b.$

By Lemma \ref{skun}, c)
the above assumptions imply that
$T_{nm/4}(a)=\pm 1$. 
Furthermore, since $$T_{mn/4}(a)=T_{m/2}(T_{n/2}(a))=T_{m/2}(-T_{n/2}(b))=-T_{mn/4}(b),$$
it follows from $W_3(a)=W_3(b)$ and 
the first equality in \eqref{lak} that $R(1)=0$. 

Changing now $n$ to $2n$ and $m$ to $2m$ we conclude that 
the fourth case
listed in Theorem \ref{os2} has place. \qed

\vskip 0.2cm
\noindent {\bf Remark.} 
Note that a solution of the fourth type appearing in Theorem \ref{os2} in general may not be obtained as a sum of only two reducible solutions. 
Consider for example the following in a sense simplest possible pair of $P$ and $Q$ as in 4)
\be \la{xruk} P=z^2(z^2-1)^2\circ T_{mn},
\ \ \ Q= T_{2m}+ T_{2n}+zR(z^2)\circ T_{mn}.
\ee 
Assume additionally that $m,n$ are different prime greater than three, and 
show that $Q$ can not be represented as a sum of two reducible solutions.

Observe first that $P$ is not linearly equivalent to a Chebyshev polynomial since $\pm 1$ are critical values of $T_{mn}$ and at the same time are critical points of the polynomial $z^2(z^2-1)^2$ implying that $P$ has critical points of multiplicity four.
Further, show that, up to the linear equivalence, compositional 
right factors of $P$ are $T_2$, $T_m,$ $T_n$, $T_{2nm},$ $T_{2n}$, $T_{2m},$ $T_{mn}$, or $zR(z^2)\circ T_{mn}$. Indeed, all the polynomials above are clearly right factors of $P$. On the other hand, it follows from Lemma \ref{c3} applied to the equality $$P=z^2(z^2-1)^2\circ T_{mn}=V\circ W$$ that if $W$ is a compositional right factor of $P$, then 
either $mn\vert \deg W$ 
or $W$ is linearly equivalent 
to a Chebyshev polynomials. 
In the first case Theorem \ref{r1} yields that 
$W=U\circ T_{mn},$ where $U$ is a right factor of $z^2(z^2-1)^2$, 
implying that 
$W$ is linearly equivalent either to $T_{2nm}$ or $zR(z^2)\circ T_{mn}$.
On the other hand, taking into account that $n$, $m$ are prime greater than three, in the second case $W$ is linearly equivalent either to one of the Chebyshev polynomials listed above either to a Chebyshev polynomial whose order is divisible by three.
However, the last case is not possible for otherwise $T_3$ also would be a right factor of $P$ and Lemma \ref{c3} 
applied to the equality
$$P=\frac{z+1}{2}\left(\frac{z-1}{2}\right)^2\circ T_{2mn}=V\circ T_3$$ would imply that  
$P$ is linearly equivalent to a Chebyshev polynomial.

The conditions imposed on $a$ and $b$ imply that among compositional right factors of $P$ only the polynomials 
$T_{2n},$ $T_{2m},$ $T_{2mn},$ and $zR(z^2)\circ T_{mn}$ satisfy the 
condition $W(a)=W(b).$ Therefore, taking into account that $T_{2mn}=T_n\circ T_{2m}=
T_m\circ T_{2n},$ we conclude that if $Q$ may be represented as a sum of two reducible solutions, then $Q$ has the form \be \la{kro} Q=V_1\circ W_{1}+V_2\circ W_{2},\ee where $W_1, W_2$ are different polynomials from the set $$S=\{T_{2n},T_{2m}, zR(z^2)\circ T_{mn}\}$$ and $V_1,$ $V_2\in \C[z].$

Furthermore, it follows from \eqref{kro} that for the polynomial 
$W_3$ from $S$ distinct from $W_1$, $W_2$ the equality  
$$W_3=(V_1-1)\circ W_{1}+(V_2-1)\circ W_{2}$$ holds.
However, the last equality is impossible since any two polynomials
$W_1,$ $W_2$ from $S$ have a common compositional right factor
which is not a compositional right factor of $W_3.$

\end{document}